\documentclass[11pt,
]{article}

\usepackage{amssymb,amsfonts,amstext,amsmath,amsthm,latexsym,mathrsfs}
\usepackage{mathbbol}

\usepackage{makeidx,mparhack,vector}

\makeindex

\usepackage{xcolor}

\input{odarx.sty}

\begin{document}

\title
{A countable definable set of reals containing no definable 
elements  
\thanks{Revised version. 
The revision includes an updated proof of Lemma 4.5 
(the density-preservation lemma for the product). }
}

\author 
%, Lyubetsky]
{
%Ali~Enayat\thanks{address, email}  
%\and
Vladimir~Kanovei\thanks{IITP RAS and MIIT,
  Moscow, Russia, \ {\tt kanovei@googlemail.com} --- contact author. 
Partial support of 
RFFI grants 13-01-00006 (of the 2014 version) and 17-01-00705 acknowledged.}  
\and
Vassily~Lyubetsky\thanks{IITP RAS,
  Moscow, Russia, \ {\tt lyubetsk@iitp.ru} 
}}

\date 
{\today}

\maketitle

\begin{abstract}
We make use of a finite support product of Jensen forcing to 
define a model in which there is a countable non-empty $\ip12$ 
set $X$ of reals containing no ordinal-definable real.\snos
{The result was strengthened in \cite{kl22}, to the effect that 
the counterexample set $X$ is a $\Eo$-equivalence class, 
or a Vitali equivalence 
class (a shift of $\mathbb Q$, the rationals), if the true reals of 
the real line $\mathbb R$ are considered.}
\end{abstract}

\parf{Introduction}
\las{cha1}

It is well-known that the existence of a non-empty $\od$ 
(ordinal-definable) set of reals $X$ with no \od\  
element is consistent with $\ZFC$; the set of all 
non-constructible reals gives an example in many generic 
models including \eg\ the Solovay model or the extension 
of $\rL$, the constructible universe, by a Cohen real.
Can such a set $X$ be countable? 

This question was initiated and briefly 
discussed at the Mathoverflow exchange 
desk in 2010\snos
{\label{snos1}
A question about ordinal definable real numbers. 
Mathoverflow, March 09, 2010. 
{\tt http://mathoverflow.net/questions/17608}. 
}
and at FOM\snos
{\label{snos2}%
Ali Enayat. Ordinal definable numbers. FOM Jul 23, 2010.
{\tt http://cs.nyu.edu/pipermail/fom/2010-July/014944.html}}
.
In particular Ali Enayat (Footnote~\ref{snos2}) conjectured that 
the problem can be solved by the   
finite-support product $\plo$ 
of countably many copies of the Jensen ``minimal $\ip12$ real 
singleton forcing'' $\dP$ defined in \cite{jenmin} 
(see also Section 28A of \cite{jechmill}).
Enayat proved that a symmetric part of the \dplo generic 
extension of $\rL$ definitely yields a model of $\zf$ 
(not a model of $\ZFC$!) 
in which there is a Dedekind-finite infinite \od\ set of 
reals with no \od\ elements. 
In fact both \dplo generic extensions and their symmetric 
submodels were considered in \cite{ena} (Theorem 3.3) with 
respect to some other questions. 

Following the mentioned conjecture, we prove 
the next theorem in this paper$:$

\bte
\lam{mt}
It is true in a\/ \dplo generic extension of\/ $\rL$, 
the constructible universe, that 
the set of\/ \dd\dP generic reals is non-empty, 
countable, and\/ $\ip12$, 
but it has no\/ $\od$ elements.
\ete

The $\ip12$ definability is definitely the best one can get 
in this context since it easily follows from the $\ip11$ 
uniformisation theorem that any non-empty $\is12$ set 
of reals definitely contains a $\id12$ element.

Jindra Zapletal\snos
{\label{snos3}%
Personal communication, Jul 31/Aug 01, 2014.} 
informed us that there is a totally different model of $\ZFC$ 
with an $\od$ \dd\Eo class $X$ containing no \od\ elements.
The construction of such a model, not yet published, 
but described to us in a brief communication, 
looks quite complicated and involves a combination of 
several forcing notions and 
some modern ideas in descriptive set theory 
recently presented in \cite{ksz}; it also does not look to 
be able to get $X$ analytically definable, let alone $\ip12$.

%\vyk{
It remains to note that a \rit{finite} \od\ set of reals 
contains only \od\ reals by obvious reasons. 
On the other hand, by a result in \cite{gl} 
there can be two \rit{sets} of reals $X,Y$  
such that the pair $\ans{X,Y}$ is \od\ but 
neither $X$ nor $Y$ is \od.  
%}

\back
The authors thank Jindra Zapletal and Ali Enayat for fruitful discussions.
\eack

\parf{Trees and perfect-tree forcing} 
\las{tre}
\label{ptf}

Let $\bse$ be the set of all strings (finite sequences) 
of numbers $0,1$.
If $t\in\bse$ and $i=0,1$ then 
$t\we k$ is the extension of $t$ by $k$. 
If $s,t\in\bse$ then $s\sq t$ means that $t$ extends $s$, while 
$s\su t$ means proper extension. 
If $s\in\bse$ then $\lh s$ is the length of $s$,  
and $2^n=\ens{s\in\bse}{\lh s=n}$ (strings of length $n$).%

A set $T\sq\bse$ is a \rit{tree} iff 
\index{tree}%
%it is an initial segment, that is, 
for any strings $s\su t$ in $\bse$, if $t\in T$ then $s\in T$.
Thus every non-empty tree $T\sq\bse$ contains the empty 
string $\La$. 
If $T\sq\bse$ is a tree and $s\in T$ then put 
$T\ret s=\ens{t\in T}{s\sq t\lor t\sq s}$. 
%this is a tree as well.

Let $\pet$ be the set of all \rit{perfect} trees 
$\pu\ne T\sq \bse$. 
\imar{pet}%
Thus a non-empty tree $T\sq\bse$ belongs to $\pet$ iff 
it has no endpoints and no isolated branches. 
Then there is a largest string $s\in T$ such that 
$T=T\ret s$; it is denoted by $s=\roo T$   
(the {\it stem\/} of $T$); 
we have $s\we 1\in T$ and $s\we 0\in T$ in this case.

Each perfect tree $T\in\pet$ defines  
$ 
[T]=\ens{a\in\dn}{\kaz n\,(a\res n\in T)}\sq\dn
$, 
the perfect set of all \rit{paths through $T$}. 
\vyk{
 then accordingly 
$T=\ctr{[T]}$, where 
$$
\ctr X=\ens{a\res n}{a\in X\land n\in\om}\sq\bse 
\quad\text{for any set}\quad
X\sq\dn.
$$
}

By a {\ubf perfect-tree forcing} we understand any set 
$\dP\sq\pet$ suct that 
\ben
\nenu
\itla{ptf1} 
$\dP$ contains the full tree $\bse$;
\imar{ptf1}

\itla{ptf2}
\label{utp} 
if $u\in T\in\dP$ then $T\ret u\in \dP$.
\imar{ptf2}
\een
Such a set $\dP$ can be considered as a forcing notion 
(if $T\sq T'$ then $T$ is a stronger condition). 
The forcing $\dP$ adds a real in $\dn$. 

%Let $\dP$ be a perfect-tree forcing notion. 
%
Let $\plo$ be the 
{\ubf product of \dd\om many copies of $\dP$ 
with finite support.} 
Thus a typical element of $\plo$ is a sequence 
$\jta=\sis{T_n}{n\in\om}$, 
where each term $T_n=\jta(n)$ belongs to $\dP$
and the set $\abs \jta=\ens{n}{T_n\ne\bse}$
(the support of $\jta$) is finite. 
We order $\plo$ componentwisely: $\jsg\leq\jta$ 
($\jsg$ is stronger) iff $\jsg(n)\sq\jta(n)$ in $\dP$ for all $n$;
$\plo$ adds an infinite sequence 
$\sis{x_n}{n<\om}$ of %mutually 
\dd\dP generic reals 
$x_n\in\dn$. 
% whose names will be $\dox_n$.

\bre
\lam{tup}
Sometimes we'll use tuples like 
$\ang{T_0,\dots,T_n}$ of trees $T_i\in\dP$ to denote 
the infinite sequence
$\ang{T_0,\dots,T_n,\bse,\bse,\bse,\dots}\in\plo$. 
\ere

\parf{Splitting construction over a perfect set forcing} 
\las{spe}

Assume that $\dP\sq\pet$ is a perfect-tree forcing notion.
The \rit{splitting construction} $\spe\dP$ over $\dP$ 
consists of all finite systems of trees of the form 
$\vpi=\sis{T_s}{s\in2^{< n}}$, where 
$n=\vys\vpi<\om$ (the height of $\vpi$) and
\ben
\nenu
\atc\atc
\itla{spe1}
each tree $T_s=\vpi(s)$ belongs to $\dP$;
\imar{spe1}

\itla{spe2}
if $s\we i\in 2^{<n}$ ($i=0,1$) then $T_{s\we i}\sq T_s$ and  
$\roo{T_{s}}\we i\sq\roo{T_{s\we i}}$ --- 
\imar{spe2}
it easily follows that 
$[T_{s\we0}]\cap [T_{s\we0}]=\pu$.
\een
The empty system $\La$ is the only one in 
$\spe\dP$ satisfying $\vys\La=0$.

Let $\vpi,\psi$ be systems in $\spe\dP$.
Say that 
\bit
\item[$-$]
$\vpi$ \rit{extends} $\psi$, symbolically $\psi\cle\vpi$, if 
$n=\vys\psi\le\vys\vpi$ and $\psi(s)=\vpi(s)$ for 
all $s\in2^{<n}$;

\item[$-$]
\rit{properly extends} $\psi$, 
symbolically $\psi\cls\vpi$, if in 
addition $\vys\psi<\vys\vpi$;

\item[$-$]
\rit{reduces} $\psi$, if $n=\vys\psi=\vys\vpi$, 
$\vpi(s)\sq\psi(s)$ for all $s\in 2^{\vys\vpi-1},$ and 
$\vpi(s)=\psi(s)$ for all $s\in 2^{<\vys\vpi-1}$.
\eit
In other words, reduction allows to shrink trees in the top 
layer of the system, but keeps intact those in the lower 
layers.

Under the above assumption \ref{utp}, there is a strictly 
\dd\cls increasing sequence $\sis{\vpi_n}{n<\om}$ in 
$\spf\dP$. 
The limit system $\vpi=\bigcup_n\vpi_n=\sis{T_s}{s\in\bse}$
then satisfies \ref{spe1} and \ref{spe2} on the whole domain 
$\bse,$ 
and in this case, $T=\bigcap_n\bigcup_{s\in2^n}T_s$ is 
still a perfect tree 
in $\pet$ 
(not necessarily in $\dP$), 
and $[T]=\bigcap_n\bigcup_{s\in2^n}[T_s]$.

Say that a tree $T$ \rit{occurs in\/ $\vpi\in\spf\dP$} if 
$T=\vpi(s)$ for some $s\in 2^{<\vys\vpi}$.

We define $\spg\dP$, 
{\ubf the finite-support product of $\spf\dP$,} 
to consist of all infinite sequences  
$\Phi=\sis{\vpi_k}{k\in\om}$, where each $\vpi_k=\Phi(k)$ 
belongs to $\spf\dP$
and the set $\abs \Phi=\ens{k}{\vpi_k\ne\La}$
(the support of $\Phi$) is finite. 

Say that a tree $T$ \rit{occurs in\/ 
$\Phi=\sis{\vpi_k}{}$} if 
it occurs in some $\vpi_k\yd k\in\abs\Phi$.

We define $\Psi\cle\Phi$ iff $\Psi(k)\cle\Phi(k)$ 
(in $\spf\dP$) for all $k$. 
Then $\Psi\cls\Phi$ means that $\Psi\cle\Phi$ and 
$\Psi(k)\cls\Phi(k)$ for at least one $k$. 
In addition we define $\Psi\cll\Phi$ iff $\abs\Psi\sq\abs\Phi$ 
and $\Psi(k)\cls\Phi(k)$ for all $k\in\abs\Phi$.

\parf{Jensen's extension of a perfect tree forcing} 
\las{jex}

Let $\zfc'$ be the subtheory of $\zfc$ including all 
axioms except for the 
power set axiom, plus the axiom saying that $\pws\om$ exists. 
(Then $\omi$ and continual sets like $\pet$ exist as well.)
Let $\cM$ be a countable transitive model of $\ZFC'$. 

Suppose that $\dP\in\cM\yd \dP\sq\pet$ is a perfect-tree  
forcing notion.  
Then the sets $\plo$, $\spf\dP$, and $\spg\dP$ belong to 
$\cM$, too. 

\bdf
\lam{dPhi}
Consider any \dd\cle increasing sequence 
$\dphi=\sis{\Phi^j}{j<\om}$ of systems 
$\Phi^j=\sis{\vpj jk}{k\in\om}\in\spg\dP$, 
\rit{generic over\/ $\cM$} in the sense that it intersects 
every set $D\in\cM\yd D\sq\spg\dP$, dense in $\spg\dP$\snos
{Meaning that for any $\Psi\in\spg\dP$ there is $\Phi\in D$ 
with $\Psi\cle\Phi$.} 
.

Then in particular it intersects every set of the form   
$$
D_k=\ens{\Phi\in\spg\dP}{\kaz k'\le k\: 
(k\le\vys{\Phi(k')}}\,. 
$$
Hence if $k<\om$ then the sequence 
$\sis{\vpj jk}{j<\om}$ of systems $\vpj jk\in\spf\dP$ is 
\rit{eventually strictly increasing}, so that 
$\vpj jk\cls \vpj {j+1}k$ 
for infinitely many indices $j$  
(and $\vpj jk=\vpj {j+1}k$ for other $j$).
Therefore there is a system of trees 
$\sis{\tx ks}{k<\om\land s\in\bse}$ in $\dP$ 
such that $\vpj jk=\sis{\tx ks}{s\in 2^{<h(j,k)}}$, 
where $h(j,k)=\vys{\vpj jk}$.
Then 
$$
\textstyle
\ufi k=\bigcap_n\bigcup_{s\in2^n}\tx ks\qand 
\uf ks\bigcap_{n\ge \lh s}
\bigcup_{t\in2^n,\:s\sq t}\tx kt
$$   
are trees in $\pet$ (not necessarily in $\dP$) 
for each $k$ and $s\in\bse;$ 
thus $\ufi k=\uf k\La$. 
In fact $\uf ks=\ufi k\cap \tx ks$ by \ref{spe2}.
\edf

\ble
\lam{uu1}
The set of trees\/ 
$\dU=\ens{\uf ks}{k<\om\land s\in\bse}$
satisfies\/ \ref{ptf2} while the union\/ 
$\dP\cup\dU$ is a perfect-tree forcing.
% of Section~\ref{ptf}.
\qed
\ele

\ble
\lam{uu2}
The set\/ $\dU$ is dense in\/ $\dU\cup\dP$. 
\ele
\bpf
Suppose that $T\in\dP$. 
The set $D(T)$ of all systems   
$\Phi=\sis{\vpi_k}{k\in\om}\in \spg\dP$,
such that $\vpi_k(\La)=T$ for some $k$, belongs to $\cM$ 
and obviously is dense in $\spg\dP$. 
It follows that $\Phi^j\in D(T)$ for some $j$, 
by the choice of $\dphi$. 
Then $\tx k\La=T$ for some $k$. 
However $\uf k\La\sq \tx k\La$. 
\epf

\ble
\lam{uu3}
If a set\/ $D\in\cM$, $D\sq\dP$ is pre-dense 
in\/ $\dP$, and\/ $U\in\dU$, then\/ $U\sqf\bigcup D$, 
that is, there is a finite\/ $D'\sq D$ with
$U\sq\bigcup D'$. 
\ele
\bpf
Suppose that $U=\uf Ks$, $K<\om$ and $s\in\bse.$ 
Consider the set $\Da\in\cM$ of all systems   
$\Phi=\sis{\vpi_k}{k\in\om}\in \spg\dP$ such that 
$K\in\abs\Phi$, $\lh s<h=\vys{\vpi_K}$, 
and for each $t\in 2^{h-1}$ there is a tree $S_t\in D$ with  
$\vpi_K(t)\sq S$.
The set $\Da$ 
is dense in $\spg\dP$ by the pre-density of $D$. 
Therefore there is an index $j$ such that $\Phi^j$ belongs 
to $\Da$.
Let this be witnessed by trees $S_t\in D\yt t\in 2^{h-1},$ 
where $\lh s<h=\vys{\vpj JK}$, so that $\vpj JK(t)\sq S_t$. 
Then 
$$
\textstyle
U=\uf Ks\sq\uf K\La\sq\bigcup_{t\in 2^{h-1}}\vpj JK(t)
\sq\bigcup_{t\in 2^{h-1}}S_t\sq\bigcup D'
$$ 
by construction, where $D'=\ens{S_t}{t\in 2^{h-1}}\sq D$ 
is finite.
\epf

\ble
\lam{uu4}
If a set\/ $D\in\cM,$ $D\sq\plo$  
is pre-dense in\/ $\plo$ then it remains pre-dense in\/ 
$\uplo$. 
\ele
\bpf\snos
{An improved argument, first appeared in 
a more complicated case in 
\cite[Theorem 6.3]{kl28}.}
Given a condition $\jta\in\uplo$, we have to prove that 
$\jta$ is compatible in $\uplo$ with a condition 
$\jsg\in D$.
For the sake of brevity, assume that $\jta=\ang{U,V}$,  
where $U=\uf ks$ and $V=\uf\ell t$ belong to $\dU$. 
The numbers $k,\ell$ can be equal or different.

Consider the set $\Da\in\cM$ of all systems   
$\Phi=\sis{\vpi_k}{k\in\om}\in \spg\dP$ such that 
$k,\ell\in\abs\Phi$ and 
there exist:
\ben
\fenu
\itla{f*}
 strings $s',t'\in\bse$ with $s\sq s'$, $t\sq t'$, 
$\lh s'<\vys{\vpi_k}$, $\lh t'<\vys{\vpi_\ell}$, 
and tuples 
$\jsg=\ang{S_0,S_1,\dots,S_{n-1}}\in\plo$,   
$\bro=\ang{R_0,R_1,\dots,R_{n-1}}\in D$
$(n\ge2)$,
such that all trees $S_i$ occur in $\Phi$, 
$S_i\sq R_i$ for all $i$, and finally 
$\vpi_k(s')=S_0$, $\vpi_\ell(t')=S_1$.
\een
The set $\Da$ 
is dense in $\spg\dP$ by the pre-density of $D$. 
Therefore there is an index $j$ such that $\Phi^j$ belongs 
to $\Da$.

Let this be witnessed by strings $s',t'\in\bse$ and 
tuples $\jsg,\jta$ as in \ref{f*}. 
By definition there exists a tuple 
$\ju=\ang{U_0,U_1,\dots,U_{n-1}}\in\dU\lom,$ 
such that $U_i\sq S_i\sq R_i$ for all $i$ --- 
hence $\ju$ is stronger than $\bro\in D$, --- and 
$U_0=\uf k{s'}$, $U_1=\uf \ell{t'}$.  
However $\uf k{s'}\sq \vpj jk(s')\cap \uf k{s}$ and 
$\uf\ell{t'}\sq \vpj j \ell (t')\cap \uf\ell{t}$ by construction. 
It follows that condition $\ju\in \dU\lom$ is stronger 
than both $\jta=\ang{U,V}$ and $\bro\in D$, as required.
\epf

\vyk{
\bpf
Given a condition $\jta\in\uplo$, we have to prove that 
$\jta$ is compatible in $\uplo$ with a condition 
$\jsg\in D$.
For the sake of brevity, assume that $\jta=\ang{U,V}$,  
where $U=\uf ks$ and $V=\uf\ell t$ belong to 
$\dU$.

Consider the set $\Da\in\cM$ of all systems   
$\Phi=\sis{\vpi_k}{k\in\om}\in \spg\dP$ such that there 
are strings $s',t'\in\bse$ with $s\sq s'$, $t\sq t'$, 
$\lh s'<\vys{\vpi_k}$, $\lh t'<\vys{\vpi_\ell}$, 
and trees $S,T\in\dP$ such that 
$\ang{S,T}\in D$ and $\vpi_k(s')\sq U\cap S$, 
$\vpi_\ell(t')\sq V\cap T$.
The set $\Da$ 
is dense in $\spg\dP$ by the pre-density of $D$. 
Therefore there is an index $j$ such that $\Phi^j$ belongs 
to $\Da$.

Let this be witnessed by $s',t'\in\bse$ and $\ang{S,T}\in D$. 
In other words, $\vpj jk(s')\sq U\cap S$ and  
$\vpj j\ell(t')\sq V\cap T$. 
However $U'=\uf k{s'}\sq \vpj jk(s')$ and 
$V'=\uf\ell{t'}\sq \vpj j \ell (t')$ by construction. 
It follows that condition $\ang{U',V'}\in \dU\lom$ is stronger 
than both $\ang{U,V}$ and $\ang{S,T}$, as required.
\epf
}

\parf{Forcing a real away of a pre-dense set} 
\las{saway}

Let $\cM$ be still a countable transitive model of $\ZFC'$ 
and $\dP\in\cM\yd \dP\sq\pet$ be a perfect-tree  
forcing notion. 
The goal of the following Theorem~\ref{K} is to prove that, 
in the conditions of Definition~\ref{dPhi}, for any 
\dd\plo name $c$ of a real in $\dn,$ it is forced by the extended 
forcing $\uplo$ that $c$ does not belong to sets $[U]$ where 
$u$ is a tree in $\dU$ --- unless $c$ is a name of one of 
generic reals $x_k$ themselves.
We begin with a suitable notation.

\bdf
\lam{rk}
A \rit{\dplo real name} is a system 
$\rc=\sis{\kc ni}{n<\om\yi i<2}$ of sets $\kc ni\sq\dP\lom$ 
such that each set $C_n=\kc n0\cup \kc n1$ is 
dense or at least pre-dense in $\plo$ 
and if $\jsg\in \kc n0$ and $\jta\in \kc n1$ then $\jsg,\jta$ are 
incompatible in $\plo$.

If a set $G\sq\plo$ is \dd\plo generic at least over 
the collection of all  sets $C_n$ then we define 
$\rc[G]\in\dn$ so that $\rc[G](n)=i$ iff $G\cap \kc ni\ne\pu$.
\edf

Thus any \dplo  real name $\rc=\sis{\kc ni}{}$ 
is a \dplo  name for a real in $\dn.$ 

Recall that $\dP\lom$ adds a generic sequence 
$\sis{x_k}{k<\om}$ of reals $x_k\in\dn$.

\bpri
\lam{proj}
Let $k<\om$. 
Define a \dplo real name $\rpi_k=\sis{\kcp nik}{n<\om\yi i<2}$
such that each set $\kcp nik$ contains a single condition 
$\kr nik\in\plo$, and $\abs{\kr nik}=\ans{k}$,  
$\kr nik(k)=\kR ni$, where 
$\kR ni=\ens{s\in\bse}{\lh s>n\imp s(n)=i}$.
Then $\rpi_k$ is a \dplo name of a real $x_k$, the $k$th 
term of a \dplo generic sequence $\sis{x_k}{k<\om}$. 
\epri

Let $\rc=\sis{\kc ni}{}$ and $\rd=\sis{\kc ni}{}$ 
be a \dplo real names. 
Say that  $\jta\in \pet\lom$:
\bit
\item
\rit{directly forces\/ $\rc(n)=i$}, 
where $n<\om$ and $i=0,1$, iff $\jta\leq\kr nik$ 
(that is, the tree $T=\jta(k)\in\pet$ satisfies 
$x(n)=i$ for all $x\in[T]$); 

\item
\rit{directly forces\/ $s\su\rc$},  
where $s\in\bse,$ iff for all $n<\lh s$, $\jta$ 
directly forces $\rc(n)=i$, where $i=s(n)$; 

\item
\rit{directly forces\/ $\rd\ne\rc$}, iff there are strings 
$s,t\in\bse,$ incomparable in $\bse$ and such that  
$\jta$ directly forces $s\su\rc$ 
and $t\su\rd$; 

\item
\rit{directly forces\/ $\rc\nin[T]$},  
where $T\in\pet$, iff there is a string $s\in\bse\bez T$ 
such that $\jta$ directly forces $s\su \rc$; 
\vyk{
\item
\,{\ubf[applicable only for $\jta\in\plo$]}
\rit{weakly\/ \dd\plo forces\/ $\rc\nin[T]$},  
iff the set of all conditions 
$\jsg\in\plo$ that directly force\/ $\rc\nin[T]$ 
is dense in $\plo$ below $\jta$.
} 
\eit

\bte
\lam{K}
In the assumptions of Definition~\ref{dPhi}, suppose that\/ 
$\rc=\break
\sis{C_m^i}{m<\om\yi i<2}\in\cM$ is a\/ \dplo real name, 
and for every\/ $k$ the set
$$
D(k)=\ens{\jta\in\dP\lom}{\jta\text{ directly forces }\rc\ne\rpi_k}
$$
is dense in\/ $\plo.$ 
Let\/ $\ju\in(\dP\cup\dU)\lom$ and\/ $U\in \dU$.
Then there is\/ 
a stronger condition\/ $\jv\in\dU\lom\yd \jv\leq\ju$, which 
directly forces\/ $\rc\nin[U]$.
\ete
\bpf
By construction $U\sq U^\dphi_k$ for some $k$; 
thus we can assume that simply $U=U^\dphi_k$. 
Let, say, $U=U^\dphi_1$.
Assume for the sake of brevity that 
$K=1$, $\abs\jta=\ans{0,1,2,3}$, and  
$\ju=\ang{U_0,U_1,U_2,U_3}\in\dU\lom$ (see Remark~\ref{tup}), 
where
$$
U_0=\uf0{t_0}\,,\quad
U_1=\uf0{t_1}\,,\quad
U_2=\uf1{t_2}\,,\quad
U_3=\uf1{t_3}\,,
$$
and  $t_0,t_1,t_2,t_3$ are strings in $\bse.$ 

There is an index $J$ such that the system 
$\Phi^J=\sis{\vpj Jk}{k\in\om}$ satisfies 
$\vys{\vpj J0}>\tmax\ans{\lh{t_0},\lh{t_1}}$ and 
$\vys{\vpj J1}>\tmax\ans{\lh{t_2},\lh{t_2}}$, 
so that the trees 
$$
T_0=\vpj J0(t_0)=\tx 0{t_0},\;
T_1=\vpj J0(t_1)=\tx 0{t_1},\;
T_2=\vpj J1(t_2)=\tx 1{t_2},\;
$$ 
and $T_3=\vpj J1(t_3)=\tx 1{t_3}$
in $\dP$ are defined and condition 
$\jta=\ang{T_0,T_1,T_2,T_3}$  
belongs to $\plo$. 
Note that $\ju\leq\jta$.  

Consider the set $\cD$ of all systems 
$\Phi
=\sis{\vpi_k}{k\in\om}
\in \spg\dP$ such that 
$\Phi^J\cle\Phi$ and there is a condition 
$\jsg=\ang{S_0,\dots,S_n}\in\plo\yt\jsg\leq\jta$
(\ie, stronger that $\jta$), such that
\ben
\nenu
\atc\atc\atc\atc
\itla{nen1}
$\jsg$ directly forces $\rc\nin[T]$, where 
$T=\bigcup_{s\in 2^{h_1-1}}\vpi_1(s)$ and 
$h_k=\vys{\vpi_k}$;

\itla{nen2}
each tree $S_i$ occurs in $\Phi$ (see Section~\ref{spe});

\itla{nen3}
more specifically, 
$S_0=\vpi_0(s_0)\yt S_1=\vpi_0(s_1)\yt S_2=\vpi_1(s_2)
\yt S_3=\vpi_1(s_3)$, 
where $s_0,s_1\in 2^{h_0-1}$, $s_2,s_3\in 2^{h_1-1}$, 
and $t_i\sq s_i$, $i=0,1,2,3$. 
\een

\ble
\lam{Ka}
$\cD$ is dense in\/ $\spg\dP$ above\/ $\Phi^J$. 
\ele
\bpf
Consider any system $\Phi=\sis{\vpi_k}{k\in\om}\in\spg\dP$ 
with $\Phi^J\cle\Phi$; the goal is to define a system 
$\Phi'\in\cD$ such that $\Phi\cle\Phi'$.
We can assume that in fact $\Phi^J\cll\Phi$; 
then any system $\Phi'\in\spg\dP$ which is a reduction of 
$\Phi$ still satisfies $\Phi^J\cll\Phi'$ and $\Phi^J\cle\Phi'$. 
Let $h_0=\vys{\vpi_0}$ and $h_1=\vys{\vpi_1}$. 
Then by the assumption $\vys{\vpj J0}<h_0$ and 
$\vys{\vpj J1}<h_1$ strictly.

Pick strings $s_0,s_1\in 2^{h_0-1}$ and $s_2,s_3\in 2^{h_1-1}$ 
satisfying $t_i\su s_i\yt i=0,1,2,3$.
Consider the condition 
$\jrho=\ang{R_0,R_1,R_2,R_3,R_4,\dots,R_N}\in\plo,$ where 
$N=1+2^{n_1}$, 
$R_0=\vpi_0(s_0)$, 
$R_1=\vpi_0(s_1)$, 
$R_2=\vpi_1(s_2)$, 
$R_3=\vpi_1(s_3)$, 
and $\ans{R_4,\dots,R_N}$ is an arbitrary enumeration of 
$\ens{\vpi_1(s)}{s\in 2^{n_1-1}\yi s\ne s_2,s_3}$. 

It follows from the density of sets $D(k)$ that there is a 
stronger condition 
$\jsg=\ang{S_0,S_1,S_2,S_3,\dots,S_N,\dots,S_M}\in\plo,$ 
where $M\ge N$ and $S_i\sq R_i$ for all $i\le N$, which 
directly forces $\rc\ne\rpi_k$ for all $k=2,\dots,N$. 
Then there exist strings $u,v_2,\dots,v_N\in\bse$ 
such that $\jsg$ directly forces each of the formulas 
\bce
$u\su\rc$, \ and also \ 
$v_2\sq\rpi_2$, \ 
$v_3\sq\rpi_3$ , \ \dots\ , \ $v_N\sq\rpi_N$,  
\ece
and $u$ is incompatible in $\bse$ with each $v_k$.

However $\jsg$ directly forces $v_k\sq\rpi_k$ iff 
$v_k\sq\roo{S_k}$. 
We conclude that $\jsg$ directly forces\/ $\rc\nin[S]$,
where $S=\bigcup_{2\le k\le M}S_k$.

Now let $\Phi'\in \spg\dP$ be defined as follows. 
We begin with $\Phi$.\vom

{\it Step 1\/}. 
Recall that $R_0=\vpi_0(s_0)$, 
$R_1=\vpi_0(s_1)$, $R_2=\vpi_1(s_2)$, $R_3=\vpi_1(s_3)$ in $\Phi$. 
Now let $\vpi'_0(s_0)=S_0$, $\vpi'_0(s_1)=S_1$, 
$\vpi'_1(s_2)=S_2$, $\vpi'_1(s_3)=S_3$.\vom

{\it Step 2\/}.  
By construction each $R_k\yd 4\le k\le M$, was equal to some 
$\vpi_1(s_k)\yd s_k\in 2^{n_1-1}\yi s_k\ne s_2,s_3$; \ we
let $\vpi'_1(t)=S_{k}$.\vom

{\it Step 3\/}. 
Each $S_k\yd N+1\le k<M$, is a tree in $\dP$. 
Let $\mu=\tmax{\abs\Phi}$ and define a system 
$\vpi'_{\mu+k} \in\spf\dP$ so that 
$\vys{\vpi'_{\mu+k}}=1$ and $\vpi'_{\mu+k}(\La)=S'_k$.\vom

After all these changes in $\Phi$, we obtain another system 
$\Phi'=\ens{\vpi'_k}{k\in\om}$ in $\spg\dP$ which is a 
reduction of $\Phi$, hence, satisfies $\Phi^J\cle\Phi'$, 
and every tree $S_k$ in the condition 
$\jsg=\ang{S_0,S_1,S_2,S_3,\dots,S_N,\dots,S_M}$ 
occurs in $\Phi'$.
Moreover $\jsg$ witnesses that $\Phi'\in\cD$, as required.
\epF{Lemma}

Come back to the proof of the theorem.
It follows from the lemma that there is an index $j\ge J$ 
such that the system $\Phi^j=\sis{\vpj jk}{k\in\om}$ 
belongs to $\cD$, and let this be witnessed by 
a condition 
$\jsg=\ang{S_0,S_1,S_2,S_3,\dots,S_n}\in\plo$ 
satisfying \ref{nen1}, \ref{nen2}, \ref{nen3}.
In particular $\jsg\leq\jta$ by \ref{nen3}. 

Finally consider a condition 
$\jv=\ang{V_0,V_1,V_2,V_3,\dots,V_n}\in\dU\lom$ 
defined so that 
$ 
V_0=\uf0{s_0}\yt
V_1=\uf0{s_1}\yt
V_2=\uf1{s_2}\yt
V_3=\uf1{s_3}$, 
and if $4\le k\le n$ then let $V_k$ be any tree in 
$\dU$ satisfying $V_k\sq S_k$ 
(Lemma~\ref{uu2}).
Recall that $t_i\sq s_i$ for $i=0,1,2,3$ by construction, 
therefore $\jv\leq\ju$. 
On the other hand, $\jv\leq\jsg$, therefore 
$\jv$ directly forces $\rc\nin[T]$ by \ref{nen1}, where 
$T=\bigcup_{s\in 2^{h-1}}\vpj j1(s)
=\bigcup_{s\in 2^{h-1}}\tx 1s$ 
and   
$h=\vys{\vpi_1}$.
And finally by definition  
$\ufi1\sq\bigcup_{s\in 2^{h-1}}\vpj j1(s)$, so 
$\jv$ directly forces $\rc\nin[U^\dphi_1]$, as required. 
\epf

\parf{Jensen's forcing} 
\las{jfor}

In this section, 
{\ubf we argue in $\rL$, the constructible universe.}
Let $\lel$ be the canonical wellordering of $\rL$.

\bdf
[in $\rL$]
%, = Section 3 in \cite{jenmin}]
\lam{uxi}
Following \cite[Section 3]{jenmin}, 
define, by induction on $\xi<\omi$, a countable set  of trees 
$\dU_\xi\sq\pet$ satisfying \ref{ptf2} of Section \ref{tre}, 
as follows.

Let $\dU_0$ consist of all clopen trees $\pu\ne S\sq\bse$, 
including $\bse$ itself.

Suppose that $0<\la<\omi$, and countable sets 
$\dU_\xi\sq\pet$ are already defined. 
Let $\cM_\xi$ be the least model $\cM$ of $\zfc'$ of the form 
$\rL_\ka\yt\ka<\omi$, 
containing $\sis{\dU_\xi}{\xi<\la}$ and such that 
$\al<\omi^\cM$ and all sets $\dU_\xi$, $\xi<\la$, 
are countable in $\cM$.

Then $\dP_\la=\bigcup_{\xi<\la}\dU_\xi$ is countable in $\cM$, too. 
Let $\sis{\Phi^j}{j<\om}$ be the $\le_\rL$-least sequence of 
systems $\Phi^j\in\spg{\dP_\la}$, \dd\cle increasing and generic 
over $\cM_\la$, and let $\dU_\la=\dU$ be defined, 
on the base of this sequence, as 
in Definition~\ref{dPhi}. 

Modulo technical details,  $\dP=\bigcup_{\xi<\omi}\dU_\xi$ 
is the Jensen forcing of \cite{jenmin}, and 
the finite-support product 
$\plo$ is the forcing we'll use to prove Theorem~\ref{mt}. 
\edf

\bpro
[in $\rL$]
\lam{uxip}
The sequence\/ $\sis{\dU_\xi}{\xi<\omi}$ belongs to\/ 
$\id{\hc}1$.\qed
\epro

\ble
[in $\rL$]
\lam{jden}
If a set\/ $D\in\cM_\xi\yt D\sq {\dP_\xi}\lom$ is 
pre-dense in\/ ${\dP_\xi}\lom$ then it remains pre-dense in\/ 
$\plo$. 
Hence if\/ $\xi<\omi$ then\/ 
${\dU_\xi}\lom$ is pre-dense in\/ $\plo$.
\ele
\bpf
By induction on $\la\yd \xi\le\la<\omi$, 
if $D$ is pre-dense in ${\dP_\la}\lom$ then it 
remains pre-dense in 
${\dP_{\la+1}}\lom={(\dP_\la\cup\dU_\la)}\lom$ 
by Lemma~\ref{uu4}. 
Limit steps are obvious. 
To prove the second part, note that 
${\dU_\xi}\lom$ is dense in ${\dP_{\xi+1}}\lom$ by Lemma~\ref{uu2},
and $\dU_\xi$ belongs to $\cM_{\xi+1}$.
\epf

\ble
[in $\rL$]
\lam{club}
If\/ $X\sq\HC=\rL_{\omi}$ then the set\/ $W_X$ of all 
ordinals\/ $\xi<\omi$ such that\/ 
$\stk{\rL_\xi}{X\cap\rL_\xi}$ is an elementary submodel of\/  
$\stk{\rL_{\omi}}{X}$ and\/ $X\cap\rL_\xi\in\cM_\xi$ 
is unbounded in\/ $\omi$.
More generally, if\/ $X_n\sq\HC$ for all\/ $n$ 
then the set\/ $W$ of all 
ordinals\/ $\xi<\omi$, such that\/ 
$\stk{\rL_\xi}{\sis{X_n\cap\rL_\xi}{n<\om}}$ 
is an elementary submodel of\/  
$\stk{\rL_{\omi}}{\sis{X_n}{n<\om}}$ 
and\/ $\sis{X_n\cap\rL_\xi}{n<\om}\in\cM_\xi$, 
is unbounded in\/ $\omi$.
\ele
\bpf
Let $\xi_0<\omi$. 
By standard arguments, there are ordinals $\xi<\la<\omi$, 
$\xi>\xi_0$, such 
that $\stk{\rL_\la}{\rL_\xi,X\cap\rL_\xi}$ is an elementary 
submodel of $\stk{\rL_{\om_2}}{\rL_{\omi},X}$.
Then $\stk{\rL_\xi}{X\cap\rL_\xi}$ is an elementary submodel of   
$\stk{\rL_{\omi}}{X}$, of course. 
Moreover, $\xi$ is uncountable in $\rL_\la$, hence 
$\rL_\la\sq\cM_\xi$. 
It follows that $X\cap\rL_\xi\in\cM_\xi$ since 
$X\cap\rL_\xi\in\rL_\la$ by construction.
The second claim does not differ much.
\epf

\bcor
[in $\rL$, = Lemma 6 in \cite{jenmin}]
\lam{ccc}
The forcing\/ $\plo$ satisfies CCC.
\ecor
\bpf
Suppose that $A\sq\plo$ is a maximal antichain. 
By Lemma~\ref{club}, there is an ordinal $\xi$ such that 
$A'=A\cap{\dP_\xi}\lom$ is a maximal antichain in ${\dP_\xi}\lom$ 
and $A'\in\cM_\xi$. 
But then $A'$ remains pre-dense, therefore, maximal, in the 
whole set $\dP$ by Lemma~\ref{jden}. 
It follows that $A=A'$ is countable.
\epf

\parf{The model} 
\las{mod}

We consider the sets $\dP,\plo\in\rL$ (Definition~\ref{uxi})   
as forcing notions over $\rL$. 
 
\ble
[= Lemma 7 in \cite{jenmin}]
\lam{mod1}
A real\/ $x\in\dn$ is\/ $\dP$-generic over\/ $\rL$ iff\/ 
$x\in Z=\bigcap_{\xi<\omi^\rL}\bigcup_{U\in\dU_\xi}[U]$. 
\ele
\bpf
All sets $\dU_\xi$ are pre-dense in $\dP$ by Lemma~\ref{jden}. 
On the other hand, if $A\sq\dP$, $A\in\rL$ is a maximal 
altichain in $\dP$, then easily $A\sq\dP_\xi$ for some 
$\xi<\omi^\rL$ by Corollary~\ref{ccc}.
But then every tree $U\in\dU_\xi$ satisfies 
$U\sqf\bigcup A$ by Lemma~\ref{uu3}, so that 
$\bigcup_{U\in \dU_\xi}[U]\sq \bigcup_{T\in A}[T]$.
\epf       

\bcor
[= Corollary 9 in \cite{jenmin}]
\lam{mod2}
In any generic extension of\/ $\rL$, the set of all reals 
in\/ $\dn$ $\dP$-generic over\/ $\rL$ is\/ $\ip\HC1$ and\/ 
$\ip12$. 
\ecor
\bpf
Use Lemma~\ref{mod1} and Proposition~\ref{uxip}.
\epf  

\bdf
\lam{gg}
From now on, let $G\sq\plo$ be a set \dd\plo generic over $\rL$.
If $k<\om$ then 
let $G_k=\ens{\jta(k)}{\jta\in G}$, so that each $G_k$ is 
\dd\dP generic over $\rL$ and $X_k=\bigcap_{T\in G_k}[T]$ is 
a singleton $X_k=\ans{x_k}$ whose only element $x_k\in\dn$ is 
a real \dd\dP generic over $\rL$.
\edf

The whole extension $\rL[G]$ is then equal to 
$\rL[\sis{x_k}{k<\om}]$, and our goal is now to prove that 
it contains no other $\dP$-generic reals.

\ble
[in the assumptions of Definition~\ref{gg}] 
\lam{only}
If\/ $x\in\rL[G]\cap\dn$ and\/ $x\nin\ens{x_k}{k<\om}$ 
then\/ $x$ is not a\/ $\dP$-generic real over\/ $\rL$.
\ele
\bpf
Otherwise there is a condition $\jta\in\plo$ and 
a \dd\plo real name $\rc=\sis{\kc ni}{n<\om,\,i=0,1}\in\rL$ 
such that $\jta$ \dplo forces that $\rc$ is \dd\dP generic 
while $\plo$ forces that $\rc\ne\rpi_k$ for all $k$. 
(Recall that $\rpi_k$ is a \dplo real name for $x_k$.)

Let $C_n=\kc n0\cup\kc n1$; this is a pre-dense set in $\plo.$ 
It follows from Lemma~\ref{club} that there is an ordinal 
$\la<\omi$ such that 
%$\jta\in{\dP_\la}\lom$, 
each set 
$C'_n=C_n\cap {\dP_\la}\lom$ is pre-dense in ${\dP_\la}\lom,$ 
and the sequence $\sis{\xc ni}{n<\om,\,i=0,1}$ belongs to 
$\cM_\la$, where $\xc ni=C'_n\cap \kc ni$ --- 
then $C'_n$ is pre-dense in $\plo,$ too, by Lemma~\ref{jden}. 
Thus we can assume that in fact $C_n=C'_n$, that is, 
$\rc\in\cM_\la$ and $\rc$ is a \dd{{\dP_\la}\lom}real name.

Further, as $\plo$ forces that $\rc\ne\rpi_k$, the set $D_k$ 
of all conditions $\jsg\in\plo$ which directly force 
$\rc\ne\rpi_k$, is dense in $\plo$ --- for every $k$. 
Therefore, still by Lemmas \ref{club}, 
% and \ref{jden}, 
we may 
assume that the same ordinal $\la$ as above satisfies the 
following: each set $D'_k=D_k\cap{\dP_\la}\lom$ is dense in 
${\dP_\la}\lom.$ 
%and belongs .$ 
%(And pre-dense in $\plo$ by Lemma~\ref{jden}).

Applying Theorem~\ref{K} with $\dP=\dP_\la$, $\dU=\dU_\la$, 
and $\dP\cup\dU=\dP_{\la+1}$, we conclude that for each 
$U\in\dU_\la$ the set $Q_U$ of all conditions 
$\jv\in {\dP_{\la+1}}\lom$ which directly force $\rc\nin[U]$, 
is dense in ${\dP_{\la+1}}\lom.$
As obviously $Q_U\in\cM_{\la+1}$, we further conclude that 
$Q_U$ is pre-dense in the whole forcing $\plo$ 
by Lemma~\ref{jden}.
This implies that $\plo$ forces 
$\rc\nin\bigcup_{U\in \dU_\la}[U]$, 
hence, forces that $\rc$ is not \dplo generic, 
by Lemma~\ref{mod1}.
But this contradicts to the choice of $\jta$.
\epf

Finally the next lemma is a usual 
property of finite-support product forcing.

\ble
[in the assumptions of Definition~\ref{gg}] 
\lam{sym}
If\/ $k<\om$ then\/ $x_k$ is not\/ $\od$ in\/ $\rL[G]$.\qed
\ele

Now, arguing in the \dplo generic model 
$\rL[G]=\rL[\sis{x_k}{k<\om}]$, 
we observe the countable set $X=\ens{x_k}{k<\om}$ 
is exactly the set of 
all \dd\dP generic reals by Lemma~\ref{only}, hence it 
belongs to $\ip12$ by Corollary~\ref {mod2}, and finally 
it contains no $\od$ elements by Lemma~\ref{sym}.

%This completes the proof of Theorem~\ref{mt}.

\qeD{Theorem~\ref{mt}}

\bibliographystyle{plain}
{\small
%\bibliography{u}

}

\end{document}